\documentclass[11pt]{article}
\usepackage{amsmath,amssymb}
\usepackage{graphicx}
\newtheorem{theorem}{Theorem}[section]   

\newtheorem{proposition}[theorem]{Proposition}

\def\qed{\hfill $\Box$\medskip}
\def\IR{{\bf R}}   
\def\IC{{\bf C}}

\def\diag{{\rm diag}\,}

\def\conv{{\rm conv}\,}
\def\bx{{\bf x}}
\def\bv{{\bf v}}
\def\la{{\langle}}
\def\ra{{\rangle}}
\begin{document}
\baselineskip 14.8pt
\title{
Stationary Probability Vectors of Higher-order Markov Chains}
\author{Chi-Kwong Li\thanks{Department of Mathematics, 
College of William and Mary, Williamsburg, VA 23187, USA. (ckli@math.wm.edu)} 
\ and \ Shixiao Zhang\thanks{Department of Mathematics, University of Hong Kong, 
Hong Kong. (praetere@gmail.com)}}
\date{}
\maketitle

\begin{abstract}
We consider the higher-order Markov Chain, and characterize the 
second order Markov chains admitting every probability distribution vector
as a stationary vector. The result is used to construct 
Markov chains of  higher-order with the same property.
We also study conditions under which the set of stationary
vectors of the Markov chain has a certain affine dimension.
\end{abstract}


\section{Introduction}

A discrete-time Markov chain is a stochastic process with a sequence of
random variables  
$$\left\{ {{X_t},t = 0,1,2 \ldots } \right\},$$ 
which takes on values in   a discrete finite state space 
$$\la n \ra = \{1, \dots, n\}$$ for a positive integer $n$, 
such that with time independent probability 
\begin{eqnarray*}
 {p_{ij}} &=& \Pr \left( {{X_{t + 1}=i}|{X_t} = j,{X_{t - 1}} = {i_{t - 1}},{X_{t - 2}} 
= {i_{t - 2}}, \dots ,{X_1} = {i_1},{X_0} = {i_0}} \right)  \\
&=&  \Pr \left( {{X_{t + 1}}=i|{X_t} = j} \right)
\end{eqnarray*}
holds for all $i,j,{i_0}, \cdots ,{i_{t - 1}}$. 
The nonnegative matrix $P = (p_{ij})_{1 \le i, j \le n}$ is the 
transition matrix of the Markov process is column
stochastic, i.e.,  $\sum_{i=1}^n p_{ij} = 1$ for $j = 1, \dots, n$.
Denote by 
\begin{equation}\label{omega}
\Omega_n =\left\{ \bx = (x_1, \dots, x_n)^t: x_1, \dots, x_n \ge 0,
 \  \sum_{i=1}^n x_i = 1\right\}
\end{equation}
the simplex of probability vectors in $\IR^n$.
A nonnegative vector $\bx\in \Omega_n$   
is a stationary probability vector (also known as the distribution) 
of a finite Markov Chain if $P\bx = \bx$.
By the Perron-Frobenius Theory (e.g., see \cite{3,Ross}) 
every  discrete-time Markov Chain has a stationary 
probability vector, and the vector is unique if  the transition matrix
is primitive, i.e., there is a positive integer $r$ such that all entries of $P^r$
are positive. The uniqueness condition is useful when one uses numerical schemes to 
determine the stationary vectors. With the uniqueness condition, any 
convergent scheme would lead to the unique stationary vector; e.g., see \cite{5}.
  
More generally, one may consider an $m$-th order Markov chain  such that 
\begin{eqnarray*}
{p_{i,{i_1}, \cdots ,{i_m}}} 
&=& \Pr \left( {{X_{t + 1}} = i|{X_t} = {i_1},{X_{t - 1}} 
= {i_2},  \dots ,{X_1} = {i_t},{X_0} = {i_{t + 1}}} \right) \\
&=&\Pr \left( {{X_{t + 1}} = i|{X_t} = {i_1}, \cdots ,{X_{t - m + 1}} = {i_m}} \right),  
\end{eqnarray*}
where $i,{i_1}, \cdots ,{i_m} \in \la n \ra$; see \cite{1,2}.
In other words, the current state of the process depends on  $m$ past states. 
Observe that  
$$\sum_{i=1}^n p_{i, i_1, \dots, i_m} = 1, \qquad 
1 \le i_1, \dots, i_m \le n.$$
When m=1, it is just the standard Markov Chain. There are many situations
that one would use the Markov Chain models.
We refer readers to the papers \cite{1,2,ling,19,20} and the references therein.
Note that $P=(p_{i,i_1, \dots, i_m})$ is an $(m+1)$-fold tensor of $\IR^n$ 
governing the transition of states in the $m$-th order Markov chain according to the 
following rule  
$$x_i(t+1) = \sum_{1 \le i_1, \dots, i_m \le n} p_{i,i_1, \dots, i_m} x_{i_1}(t) \cdots x_{i_m}(t),
\quad i = 1, \dots, n.$$
We will call $P$ the transition probability tensor of the Markov 
chain.\footnote{As pointed out by the referee, instead of
the tensor properties of $P$, we are actually studying the 
hypermatrix of the tensor $P$ with respect to a special 
choice of basis of $\otimes ^{m+1}\IC^n$.}
A nonnegative vector $\bx = (x_1, \dots, x_n)^t \in \IR^n$ with entries summing up to 1 is a stationary  (probability distribution) vector if 
\begin{equation}\label{eq1}
x_i = \sum_{1 \le i_1, \dots, i_m \le n} p_{i_1, \dots, i_m} x_{i_1} \cdots x_{i_m},
\quad i = 1, \dots, n.
\end{equation}
By a weaker version of the Perron-Frobenius Theorem for tensors
in \cite{L} (see also \cite{CPZ,FGH}),
a stationary vector for a higher-order Markov chain always 
exists. Moreover, the stationary vector will have positive entries
if the transition tensor $P = (p_{i,i_1, \dots, i_m})$ is irreducible, 
i.e.,  there is no non-empty proper index subset 
$I \subset \left\{ {1,2, \cdots ,n} \right\}$ 
such that ${p_{i,{i_1},\dots, {i_m}}} = 0$ 
for all ${i} \in I$, and ${i_1},\dots, {i_m} \notin I$. 

Researchers have derived sufficient conditions for the
stationary vector to be unique, and proposed some iterative methods
to find the stationary vector; see \cite{CPZ,FGH,ling,L}.
In this paper, we consider an extreme situation of the problem, namely, 
every probability vector in the simplex $\Omega_n$
is a stationary vector of a higher-order Markov chain. 
In the standard (first-order) Markov chain,
this can happen if and only if $P$ is the identity matrix. 
We show that such a phenomenon may occur for a large family of higher-ordered 
Markov chains.  In particular, we characterize those 
second order Markov chains with this property.  
The result is used to study higher-order Markov chains with a similar property.\footnote{
As pointed out by the referee, this problem is 
related the Inverse Perron-Frobenius Problem: Given a distribution what are the Markov 
chains having it as a stationary distribution?
For example, one may see \cite{GB}.}

In our discussion, we always let 
$${\cal E} = \{e_1, \dots, e_n\}$$ 
denote the standard basis for $\IR^n$.
Then $\Omega_n$ is the convex hull of the set ${\cal E}$, denoted by $\conv {\cal E}$.
For any $k \in \{1, \dots, n\}$, a subset of $\Omega_n$ obtained by 
taking the convex hull of $k$ vectors from the set ${\cal E}$ 
is a face of the simplex $\Omega_n$ of affine dimension $k-1$.
We also consider higher-order Markov chains with a $(k-1)$-dimension
face of $\Omega_n$ as the set of stationary vectors.
Other geometrical features and problems concerning the set of 
stationary vectors of  higher-order Markov chains will also be mentioned.

\baselineskip 17.5pt
\section{Second Order Markov Chains}

In the following, we  characterize those second order Markov chains
so that every vector in $\Omega_n$ is a stationary vector.
Note that for a second order Markov chains 
the conditions for the stationary vector $\bx = (x_1, \dots, x_n)^t$
in (\ref{eq1}) can be rewritten as
\begin{equation} \label{eq2}
 \bx =  (x_1 P_1 + \cdots x_n P_n) \bx, 
 \end{equation}
where  for $i = 1, \dots, n$, 
\begin{equation} \label{eq3}
P_i = (p_{ris})_{1 \le r, s \le n}
\end{equation}
is a column stochastic matrix, i.e., a nonnegative matrix 
so that the sum of entries of each column is 1.
We have the following theorem.

\begin{theorem} \label{main} Suppose $P = (p_{i,i_1,i_2})$ is the transition tensor
of a second order Markov chain. Then every vector in the set $\Omega_n$
is a stationary vector  if and only if
there are nonnegative vectors $\bv_1, \dots, \bv_n \in \IR^n$ 
with entries in $[0,1]$ such that for $i = 1, \dots, n$, 
$\bv_i = (v_{i1}, \dots, v_{in})^t$ with $v_{ii} = 0$, and
$$P_i = I_n - \diag(v_{i1}, \dots, v_{in}) + e_i \bv_i^t =
{\footnotesize\left(\begin{array}{ccccccc}
 1 - v_{i1}  & & & & & &  \\
{}& \ddots  & & & & & \\
{}& & 1 - v_{i,i-1} & & & &  \\
v_{i1} & \cdots  & v_{i,i-1} & 1 & v_{i,i + 1} & \cdots & v_{in} \\
{}&{}&{}&{}&1 - v_{i,i + 1}& {} & {} \\
{}&{}&{}&{}&{}& \ddots &{}\\
{}&{}&{}&{}&{}&{}&1 - v_{in}\\
\end{array}\right)}.$$
\end{theorem}

To prove Theorem \ref{main}, we need the following detailed analysis for the
second order Markov chain when $n = 2$.

\begin{proposition} \label{n=2case}
Let
${a_1},{a_2},{b_1},{b_2} \in \left[ {0,1} \right]$.
Consider the following equation with unknown $x\in [0,1]$:
$$ \left[ x\left( {\begin{array}{*{20}{c}}{{a_1}}&{{b_1}}\\
{1 - {a_1}}&{1 - {b_1}}\end{array}} \right) + 
\left( {1 - x}\right) \left( {\begin{array}{*{20}{c}}{{a_2}}&{{b_2}}\\
{1 - {a_2}}&{1 - {b_2}}\end{array}} \right)\right]
\left( {\begin{array}{*{20}{c}}x\\{1 - x}\end{array}} \right)
=  \left( {\begin{array}{*{20}{c}}x\\1-x
\end{array}} \right).
$$
Then one of the following holds for the above equation.
\begin{itemize}
\item[{\rm (1)}] 
If ${a_1} = 1,{a_2} + {b_1} = 1,{b_2} = 0$, then every $x \in [0,1]$ is a solution.
\item[{\rm (2)}] If ${a_2} + {b_1} < 1 = a_1$, then there are two solutions in $[0,1]$,
namely, $x=1$ and $x = \frac{{{b_2}}}{{{b_2} + 1 - {a_2} - {b_1}}}$. 
\item[{\rm (3)}] If $a_2+b_1 - a_1 > a_2+b_1 -1 \ge 0 = b_2$,
then there are two solutions in $[0,1]$,
namely, $x = 0$ and $x = \frac{a_2+b_1-1}{a_2+b_1-a_1}$.
\item[{\rm (4)}] Otherwise, there is a unique solution in $[0,1]$ determined as follows.

If ${a_1} - {a_2} - {b_1} + {b_2} = 0$, 
then $x = \frac{{{b_2}}}{{2{b_2} + 1 - {a_2} - {b_1}}}$.

If ${a_1} - {a_2} - {b_1} + {b_2} \ne 0$, 
then 
\[x= \frac{{2{b_2} + 1 - {a_2} - {b_1} - \sqrt \Delta}}
{{2\left( {{a_1} - {a_2} - {b_1} + {b_2}} \right)}}\] 
with $\Delta  = {\left( {2{b_2} + 1 - {a_2} - {b_1}} \right)^2} - 
4{b_2}\left( {{a_1} - {a_2} - {b_1} + {b_2}} \right) = {\left( {1 - {a_2} - {b_1}} \right)^2} + 4{b_2}
\left( {1 - {a_1}} \right) \ge 0$.
\end{itemize}
\end{proposition}

\it Proof. \rm  Let 
$$f(x) = \left( {{a_1} - {a_2} - {b_1} + {b_2}} \right){x^2} 
+ \left( {{a_2} + {b_1} - 2{b_2}} \right)x + {b_2}$$
be the first entry of the vector
$$ \left[ x\left( {\begin{array}{*{20}{c}}{{a_1}}&{{b_1}}\\
{1 - {a_1}}&{1 - {b_1}}\end{array}} \right) + 
\left( {1 - x}\right) \left( {\begin{array}{*{20}{c}}{{a_2}}&{{b_2}}\\
{1 - {a_2}}&{1 - {b_2}}\end{array}} \right)\right]
\left( {\begin{array}{*{20}{c}}x\\{1 - x}\end{array}} \right).
$$
We need only solve $f\left( x \right) = x$ with $x \in [0,1]$.
Then the equation corresponding to the second entry will also satisfy. Set
\[g(x) = f(x)-x = \left( {{a_1} - {a_2} - {b_1} + {b_2}} \right){x^2} 
+ \left( {{a_2} + {b_1} - 2{b_2} - 1} \right)x + {b_2} = 0.\]
Then $g\left( 0 \right) = {b_2} \ge 0$ and $g\left( 1 \right) = {a_1} - 1 \le 0$. 
By the Intermediate Value Theorem, there is at least one ${x_0} \in \left[ {0,1} \right]$ 
such that $g\left( {{x_0}} \right) = 0$.
Let 
$$\Delta  = {\left( {2{b_2} + 1 - {a_2} - {b_1}} \right)^2} 
- 4{b_2}\left( {{a_1} - {a_2} - {b_1} + {b_2}} \right) 
= {\left( {1 - {a_2} - {b_1}} \right)^2} + 4{b_2}\left( {1 - {a_1}} \right) \ge 0.$$

Suppose  ${a_1} - {a_2} - {b_1} + {b_2} = 0$. The quadratic equation reduces to 
$\left( {{a_2} + {b_1} - 2{b_2} - 1} \right)x + {b_2} = 0$. 
If ${{a_2} + {b_1} - 2{b_2} - 1 = 0}$, then one can readily check that condition (1) holds.
If ${{a_2} + {b_1} - 2{b_2} - 1 \ne 0}$, then the first case of condition (4) holds.

Suppose ${a_1} - {a_2} - {b_1} + {b_2} \ne 0$. 
If $g(0) > 0$ and $g(1) < 0$, then the quadratic function $g(x)$ can 
only have one solution in $[0,1]$. 
If ${{a_1} - {a_2} - {b_1} + {b_2}} > 0$, then $g(x) \rightarrow \infty$ as 
$x \rightarrow \infty$. Since $g(1) < 0$, the larger root of $g(x) = 0$ equals
$\frac{{2{b_2} + 1 - {a_2} - {b_1} + \sqrt \Delta}}
{{2\left( {{a_1} - {a_2} - {b_1} + {b_2}} \right)}}$
will be larger than 1.
Hence, the second case of condition (4) holds.
If ${{a_1} - {a_2} - {b_1} + {b_2}} < 0$, then $g(x) \rightarrow -\infty$ as 
$x \rightarrow -\infty$. Since $g(0) > 0$,  the smaller root of $g(x) = 0$ equals
$\frac{{2{b_2} + 1 - {a_2} - {b_1} + \sqrt \Delta}}
{{2\left( {{a_1} - {a_2} - {b_1} + {b_2}} \right)}}$
will be smaller than 0. 
Hence, the second case of condition (4) holds.

Suppose $0 = g(1) = a_1 - 1$. Then $g(x)$ will have another solution in $[0,1]$
if and only if $a_1 - a_2 - b_1+b_2 = 1-a_2 - b_1 + b_2 \ge 0$. This happens
if and only if condition (2) holds.

Suppose $0 = g(0) = b_2$ and $0\ne g(1)$. 
Then $g(x)$ will have another solution in $[0,1]$
if and only if $a_1-a_2-b_1+b_2 = a_1 - a_2 -b_1 < 0$ and the 
maximum of $g(x)$ is attained at a positive number $x$. This happens if and only
if condition (3) holds.
\qed

\medskip\noindent
{\bf Proof of Theorem \ref{main}}.
The sufficiency can be readily checked.
We focus on the necessity.
Note that Proposition \ref{n=2case} covers the case when $n = 2$.
We will use an inductive argument. It is illustrative to see the case when $n = 3$.
Consider the system
\[
\left[{x_1}\left( {\begin{array}{*{20}{c}}{{a_{11}}}&{{a_{12}}}&{{a_{13}}}\\
{{a_{21}}}&{{a_{22}}}&{{a_{23}}}\\{{a_{31}}}&{{a_{32}}}&{{a_{33}}} \end{array}} \right) + 
{x_2}\left( {\begin{array}{*{20}{c}}{{b_{11}}}&{{b_{12}}}&{{b_{13}}}\\
{{b_{21}}}&{{b_{22}}}&{{b_{23}}}\\{{b_{31}}}&{{b_{32}}}&{{b_{33}}}
\end{array}} \right)
+ {x_1}\left( {\begin{array}{*{20}{c}}{{c_{11}}}&{{c_{12}}}&{{c_{13}}}\\
{{c_{21}}}&{{c_{22}}}&{{c_{23}}}\\{{c_{31}}}&{{c_{32}}}&{{c_{33}}}\end{array}} \right)\right]
\bx = \bx.\]
If we set the third entry of the stationary vector $\bx$ to be 0, 
then we can have infinitely many solutions of the form
$\bx = \left( {\begin{array}{*{20}{c}}x\\{1 - x}\\0\end{array}} \right)$ 
with $x \in \left[ {0,1} \right]$. By the 2-by-2 case,  this happens 
if and only if the sub-matrices 
$\left( {\begin{array}{*{20}{c}}{{a_{11}}}&{{a_{12}}}\\{{a_{21}}}&{{a_{22}}}\end{array}} \right)$ 
and $\left( {\begin{array}{*{20}{c}}{{b_{11}}}&{{b_{12}}}\\{{b_{21}}}&{{b_{22}}}\end{array}} \right)$ 
are of the form  
$\left( {\begin{array}{*{20}{c}}1&{{a_{12}}}\\0&{1 - {a_{12}}}\end{array}} \right),
\left( {\begin{array}{*{20}{c}}{1 - {a_{12}}}&0\\{{a_{12}}}&1\end{array}} \right)$.
Similarly, setting the second entry of $\bx$ to be 0,  we see that  the submatrices 
$\left( {\begin{array}{*{20}{c}}{{a_{11}}}&{{a_{13}}}\\{{a_{31}}}&{{a_{33}}}\end{array}} \right)$ 
and $\left( {\begin{array}{*{20}{c}}{{c_{11}}}&{{c_{13}}}\\{{c_{31}}}&{{c_{33}}}\end{array}} \right)$
are of the form $\left( {\begin{array}{*{20}{c}}1&{{a_{13}}}\\0&{1 - {a_{13}}}\end{array}} \right),
\left( {\begin{array}{*{20}{c}}{1 - {a_{13}}}&0\\{{a_{13}}}&1\end{array}} \right)$.
Finally,  setting the first entry of $\bx$ to be 0, we see that the sub-matrices 
$\left( {\begin{array}{*{20}{c}}{{a_{22}}}&{{a_{23}}}\\{{a_{32}}}&{{a_{33}}}\end{array}} \right)$ 
and $\left( {\begin{array}{*{20}{c}}{{c_{22}}}&{{c_{23}}}\\{{c_{32}}}&{{c_{33}}}\end{array}} \right)$
are of the form  $\left( {\begin{array}{*{20}{c}}1&{{a_{23}}}\\0&{1 - {a_{23}}}\end{array}} \right),
\left( {\begin{array}{*{20}{c}}{1 - {a_{23}}}&0\\{{a_{23}}}&1\end{array}} \right)$. 
Thus, the three matrices in the equation are of the form
 $$\left( {\begin{array}{*{20}{c}}
1&{{a_{12}}}&{{a_{13}}}\\ 0&{1 - {a_{12}}}&0\\ 0&0&{1 - {a_{13}}} \end{array}} \right), \  
\left( {\begin{array}{*{20}{c}} {1 - {a_{12}}}&0&0\\ {{a_{12}}}&1&{{a_{23}}}\\ 
0&0&{1 - {a_{23}}} \end{array}} \right),
\  \left( \begin{array}{*{20}{c}} {1 - {a_{13}}}&0&0\\ 0&{1 - {a_{23}}}&0\\ 
a_{13}&a_{23}& 1 \end{array} \right).$$

More generally, suppose the result holds for the $(n-1)$-dimension case.
Consider the $n$-dimension case, and the equation 
$$(x_1 P_1 + \cdots + x_nP_n)\bx = \bx \qquad \hbox{ with } \bx \in \Omega_n.$$
Let $j \in \{1, \dots, n\}$. Setting the $j$-th entry of 
$\bx = (x_1, \dots, x_n)^t$ to be zero,   we see that for $i \ne j$, the $(n-1)$ sub-matrix 
of $P_i$ obtained by deleting its $j$th row
and $j$th column has the form
$$I_{n-1} - \diag(a_{i,1}, \dots, a_{i,j-1}, a_{i,j+1}, \dots, a_{i,n})
+ \hat e_i (a_{i,1}, \dots, a_{i,j-1}, a_{i,j+1}, \dots, a_{i,n}),$$
where $\hat e_i$ is obtained from $e_i$ by removing the $j$th entry for 
$i = 1, \dots, n$.
Combining the information for different $j = 1, \dots, n$, and 
$i = 1, \dots, j-1, j+1, \dots,n$, we see that the matrices ${P_1}, \cdots ,{P_n}$ 
have the asserted form. \qed

Theorem \ref{main} shows that
it is possible for a second order Markov chain to have many stationary vectors.
In  previous study \cite{CPZ,FGH,ling,L}, researchers obtained sufficient conditions for 
a higher-order Markov chain to have a unique stationary vector.
Here we construct a family of examples of second-order Markov chains
such that one of the following holds.

(a) There are exactly $k$ stationary vectors for a given 
$k \in \{1,\dots, n+1\}$.
 
(b) The set of stationary vectors is a $k$ dimensional face
of $\Omega_n$ for $k = 1, \dots, n-2$.

(c)  The set of stationary vectors is a disconnected set 
equal to the union of a $k$ dimensional face
of $\Omega_n$ and $\{(\sum_{j=1}^n e_j)/n\}$,
for $k = 1, \dots, n-2$.

\begin{theorem} \label{main2}
Suppose $n > 2$ and a second order Markov chain with transition tensor 
$P = (p_{i,i_1,i_2})$. Let $P_i = (p_{ris})_{1 \le r, s \le n}$ for $i = 1, \dots, n$.
Let $k \in \{1, \dots, n\}$ and $f_k = (e_1 + \cdots + e_k)/k$. 
If every column of $P_i$ equals $f_k$, then 
$f_k$ is the only stationary vector of the Markov chain.

\begin{itemize} 
\item[{\rm (1)}] 
If $k = 2$, replace the first column of  $P_1$ by $e_1$ and all the columns of $P_2$ by
$e_2$.
Then the resulting Markov chain has 2 stationary vectors, namely, 
$e_1$ and $e_2$.

\item[{\rm (2)}]
If $2 < k \le n$, replace the $i$th column of $P_i$ by $e_i$ and all other columns by
$e_k$ for $i = 1, \dots, k$. Then the resulting
Markov chain has $k$ stationary vectors, namely, 
$e_1, \dots, e_k$. 

\item[{\rm (3)}] 
Suppose $k = n$. If we replace the $i$th column of $P_i$ by $e_i$ for all 
$i = 1, \dots, n$,  then 
the resulting Markov chain has $n+1$ stationary vectors, namely, 
$e_1, \dots, e_n$ and $f_n$.  

\item[{\rm (4)}] If $k \in \{2, \dots, n-1\}$ and we replace the first $k$ 
columns of $P_i$ by $e_i$ for $i = 1, \dots, k$, 
then the set of stationary vectors for 
the Markov chain equals  
$\conv \{e_1, \dots, e_k\}$.

\item[{\rm (5)}] Suppose $k \in \{2, \dots, n-1\}$ and we reset the matrices
$P_1, \dots, P_n$ so that
the first $k$ columns of $P_i$ equal 
$$v_i = \begin{cases}  e_i & \hbox{ if }\ $i = 1, \dots, k$,\\
(e_{k+1} + \cdots + e_n)/(n-k) & \hbox{ if }\ $i = k+1, \dots, n$,\end{cases}$$ 
and all other columns equal to $f_n$.
Then the set of stationary vectors for 
the Markov chain equals  
$\{f_n\}\cup \conv \{e_1, \dots, e_k\}$.
\end{itemize}
\end{theorem}
 
\it Proof. \rm 
Suppose $k \in \{1, \dots, n\}$ and  
every column of $P_i$ equals $f_k$.
Then $\bx=(x_1, \dots, x_n)^t \in \Omega_n$ satisfies
$$\bx = (x_1P_1 + \cdots + x_nP_n)\bx = (x_1+\cdots + x_k) f_k$$
if and only if $x_{k+1} = \cdots = x_n = 0$ and $x_1 = \cdots = x_k = 1/k$.

\medskip
(1) Suppose $k = 2$, and we replace $P_1$ and $P_2$ as suggested.
Then $\bx \in \Omega_n$ satisfies 
$$\bx = (x_1 P_1 + \cdots + x_nP_n)\bx$$
if and only if 
$x_3 = \cdots = x_n = 0$ and 
$$x_1 \begin{pmatrix} 1 & 1/2 \cr 0 & 1/2\cr\end{pmatrix}+
x_2 \begin{pmatrix} 0 & 0 \cr 1 & 1\cr\end{pmatrix} = 
\begin{pmatrix} x_1 \cr x_2\cr\end{pmatrix}.$$
By Proposition \ref{n=2case}, $x_1 = 1$ or $x_2 = 1$.
So, the Markov chain has two stationary vectors $e_1$ and $e_2$.

\medskip
(2) Suppose $k > 2$, and 
the $i$th column of $P_i$ is replaced by $e_i$ 
and replace all other columns by $e_k$ for $i = 1, \dots, k$.
Direct checking shows that $e_1, \dots, e_k$ and $f_k$ are stationary vectors
of the Markov chain. Conversely,
suppose $\bx=(x_1, \dots, x_n)^t \in \Omega_n$ satisfies
$$\bx = (x_1P_1 + \cdots + x_nP_n)\bx.$$
Then $x_{k+1} = \cdots = x_n = 0$,
$$x_k = \sum_{j=1}^{k-1} x_j(1-x_j) + x_k,  \quad \hbox{ and } \quad  x_j = x_j^2 \quad 
\hbox{ for } j = 1, \dots, k-1.$$
Thus, $x_j \in \{0, 1\}$ so that $\bx = e_j$ if  $x_j = 1$
for any $j=1,\dots,k-1$.
If $x_1 = \cdots = x_{k-1} = 0$, then $x_k$ is the only nonzero entry and $\bx = e_k$.

\medskip
(3) Suppose $k = n$ and  we replace the $i$th column of $P_i$ by $e_i$ and all
$i = 1, \dots, n$.
Direct computation shows that $e_1, \dots, e_n$ and $f_n$ are stationary vectors.
Conversely, suppose $\bx = (x_1, \dots, x_n)^t \in \Omega_n$ satisfies
$$\bx = (x_1P_1 + \cdots + x_n P_n)\bx.$$
Then  
$$x_i  = \frac{1}{n} \left(\sum_{1 \le i,j\le n} x_ix_j - \sum_{j=1}^k x_j^2\right) + x_i^2 
= \frac{1}{n} \left(1 - \sum_{j=1}^k x_j^2\right) + x_i^2.$$
Let $\ell = \frac{1}{n} \left(1- \sum_{j=1}^k x_j^2\right)$ and 
consider two cases.

{\bf Case 1.} If $\ell = 0$, then $x_i \in \{0,1\}$ for each $i = 1, \dots, n$.
Thus, we have $\bx \in \{e_1, \dots, e_n\}$.

{\bf Case 2.} Suppose $\ell > 0$. Because $x_i^2 - x_i + \ell = 0$, we see that
 $x_i = \left(1\pm \sqrt{1-4\ell}\right)/2$.
If at least one of the $x_i$'s equals 
$\left(1+ \sqrt{1-4\ell}\right)/2$, then by the fact that $n > 2$,   
$$
1 =  \sum_{j=1}^n x_j  
\ge \frac{ \left(1+\sqrt{1-4\ell}\right)}{2} + (n-1)\frac{\left(1- \sqrt{1-4\ell}\right)}{2}
= 1 + (n-2) \frac{\left(1- \sqrt{1-4\ell}\right)}{2} > 1,
$$
which is a contradiction. Thus,
$x_i = \left(1- \sqrt{1-4\ell}\right)/2$ for each $i = 1,\dots, n$, and hence
$\bx = f_n$.

\medskip
(4) Clearly, every vector 
in $\conv \{e_1, \dots, e_k\}$ is a stationary vector of the Markov chain.
Conversely, suppose $\bx = (x_1, \dots, x_n)^t \in \Omega_n$ is a stationary vector.
Then
$$\bx = (x_1 P_1 + \cdots + x_n P_n) \bx$$
implies that $x_{k+1} = \cdots = x_n = 0$, and $x_1, \dots, x_k$ can be any nonnegative
numbers summing up to one.

\medskip
(5) One readily checks that $f_n$ and every vector 
in $\conv \{e_1, \dots, e_k\}$  is a stationary vector of the Markov chain.
Conversely, suppose $\bx = (x_1, \dots, x_n)^t \in \Omega_n$ is a stationary vector.
If  $\beta = (\sum_{j=k+1}^n x_j)$, then
\begin{eqnarray*}
\bx &=& (x_1 P_1 + \cdots + x_n P_n) \bx \\
&=& (1-\beta)(x_1e_1 + \cdots + x_k e_k) + \beta f_n
+ (1-\beta)\beta(e_{k+1} + \cdots + e_n)/(n-k).
\end{eqnarray*}
Thus, $x_{k+1} = \cdots = x_n$.
If all of them are zero, then $x_1, \dots, x_k$ can be any nonnegative 
numbers summing up to one. If $x_{k+1} = \cdots = x_n = r > 0$, then 
$x_i = (1-\beta)x_i + \beta/n$ so that $x_i = 1/n$ for  $i = 1, \dots, k$.
If follows that $x_{k+1} = \cdots = x_n = r = 1/n$ also.
\qed

Next, we obtain a result illustrating some additional geometrical feature of 
the set of stationary vectors of a second order Markov chain.

\begin{proposition} \label{1-face}
Consider the following equation for the stationary vectors of a second order Markov chain:
$$(x_1 P_1 + \cdots + x_n P_n)\bx = \bx.$$
Suppose the Markov chain has two stationary vector of the form
$xe_i + (1-x)e_j$ and $y e_i + (1-y) e_j$ for some $x, y \in (0,1)$ and 
$1 \le i < j \le n$. Then every vector of the form $z e_i + (1-z) e_j$
with $z \in [0,1]$ is a stationary vector of the Markov chain.
\end{proposition}

\it Proof. \rm If $n = 2$, the result follows from Proposition \ref{n=2case}.
Suppose $n \ge 3$. The hypothesis of the proposition implies that the
2-by-2 submatrices of $P_i$ and $P_j$ lying in rows and columns indexed by
$i$ and $j$ have the form
 $\left( {\begin{array}{*{20}{c}}1&a\\0&{1 - a}\end{array}} \right)$ 
and $\left( {\begin{array}{*{20}{c}}{1 - a}&0\\a&1\end{array}} \right)$.
It follows that every vector of the form $z e_i + (1-z) e_j$
with $z \in [0,1]$ is a stationary vector of the Markov chain. \qed

Proposition \ref{1-face} asserts that 
if the set of stationary vectors of a second order Markov chain
contains two interior points of a 1-dimensional face of $\Omega_n$,
then every vector in the 1-dimensional face is a stationary vector.

We conjecture that if the set of stationary vectors of a second order Markov
chain contains $k$ interior points of a $(k-1)$ dimensional face of the simplex
$\Omega_n$, then every vector in the $(k-1)$ dimensional face is a stationary  
vector.

\section{Higher-Order Markov Chains}

In this section, we use the results in Section
2 to construct higher-order Markov chains so that

(I) every vector in $\Omega_n$ is a stationary vector, and 

(II) the set of stationary vectors have different affine dimensions.

\medskip\noindent
We will identify a transition probability tensor
$P = (p_{i,i_1, \dots, i_m})$ as the $n\times n^m$ hypermatrix 
with row index $i = 1, \dots, n$, and column indexes 
$i_1 \cdots i_m$ with $i_1, \dots, i_n \in\la n\ra=\{1, \dots, n\}$
arranged in lexicographic order. For example, for $n = 2$ and $m = 3$, 
the row indexes are $1,2$, and the column indexes are 
$111, 112, 121, 122, 211, 212, 221, 222$. 
We will use the tensor (Kronecker) product notation for vectors in $\Omega_n$.
For example, 
$$\bx^{(3)} = \bx \otimes \bx \otimes \bx = (x_1x_1x_1, x_1x_1x_2, x_1x_2 x_1, x_1x_2x_2,
x_2x_1x_1, x_2x_1x_2, x_2x_2 x_1, x_2x_2x_2)^t .$$
The the stationary vector condition can be represented as the following matrix equation:
$$P \bx^{(m)} = \bx.$$
As pointed out by the referee, the aobove displayed equation 
is precisely the definition of an $L^2$-eigenpair in \cite{L} or equivalently 
a $Z$-eigenpair in \cite{Q}.

\medskip
We first consider Markov chains satisfying condition (I).  
We will illustrate the construction for the third order 
Markov chains for $n = 2$, and then describe the general construction.

\medskip\noindent
{\bf First order Markov chains.} Every vector in  $\Omega_2$ is a stationary 
vector if and only if $P = I_2$.

\medskip\noindent
{\bf Second order Markov chains.}  
We can use two copies the first order chain
$I_2$ 
to produce $\tilde P = [I_2 | I_2]$ so that
$$\tilde P \bx^{(2)} = [I_2 |I_2] \bx^{(2} =  \bx
\quad \hbox{ with } \bx = (x_1x_1, x_1x_2, x_2x_1, x_2x_2)^t.$$
Observe that the second and third entries on $\bx^{(2)}$ are the same, so one can permute
the second the third columns of $\tilde P = [I_2,I_2]$ to get 
$\tilde P_1 = \begin{pmatrix}1 & 1 & 0 & 0 \cr 0 & 0 & 1 & 1 \cr\end{pmatrix}$ so that
$\tilde P_1 \bx^{(2)} = \bx$ for every $\bx \in \Omega_2$.
Evidently, for every $a \in [0,1]$, $\tilde P_a = a\tilde P + (1-a)\tilde P_1$ will satisfy
$\tilde P_a \bx^{(2)} = \bx$. In fact, we have shown that these are all possible transition
tensors have the desired property.

\medskip
\noindent
{\bf 
Third order Markov chains.} Suppose $\tilde P = [Q|Q]$, where 
$Q = \begin{pmatrix}q_{111} & q_{112} & q_{121} & q_{122} \cr 
q_{211} & q_{212} & q_{221} & q_{222}\cr\end{pmatrix}$
satisfies $Q \bx^{(2)} = \bx$ for every $\bx \in \Omega_2$.
Then $P\bx^{(3)} = \bx$ for every $\bx \in \Omega_2$.
Now, observe that the entries of $\bx^{(3)}$ indexed by $112,121,211$ are all equal to 
$x_1^2 x_2$. So, we can permute the columns of $\tilde P$ indexed by $112,121,211$ in  
$6 (=3!)$ different ways to get matrices $\tilde P_1$ satisfying $\tilde P_1 \bx^{(3)} = \bx$.
Similarly, we can permute the columns of $\tilde P$ indexed by $122,212,221$ in 
6 different ways to get matriices $\tilde P_2$ satisfying
$\tilde P_2 \bx^{(3)} = \bx$. As a result,
we get $6^2 = 36$ matrices with the desired property.
Now, we can take convex combination these matrices to get a 
large family of matrices with the desired property.

\medskip
One easily extends the above idea to obtain the following.

\begin{theorem}
Suppose $m,n \ge 2$, and $P$ is a transition probability
tensor $P$ represented as an $n \times n^{m-1}$ matrix such that
$P \bx^{(m-1)} = \bx$ for all $\bx \in \Omega_n$.
Let $\tilde P = [P | \cdots | P] = {\bf 1} \otimes P$ with ${\bf 1} = (1, \dots, 1)\in 
\IR^{1\times n}$.
Then 
\begin{equation}\label{eqa}
\tilde P \bx^{(m)} = ({\bf 1} \otimes P)(\bx\otimes \bx^{(m-1)}) 
= ({\bf 1} \bx)\otimes (P \bx^{(m-1)}) = \bx \qquad \hbox{ for all } \bx \in \Omega_n.
\end{equation}
Moreover, one can permute the columns of $\tilde P$ 
corresponding to the entries in the vector $\bx^{(m)}$ with the same values:
$x_1^{m_1} \cdots x_n^{m_n}$ for all nonnegative sequence $(m_1, \dots, m_n)$ 
with $m_1 + \cdots + m_n = m$ to yield other Markov chains
satisfying {\rm (\ref{eqa})}; in addition, taking convex combination of
these matrices will also result in Markov chains satisfying {\rm (\ref{eqa})}.
\end{theorem}

Note that there are  
$${m \choose m_1, m_2, \dots, m_n} = \frac{m!} {m_1! \cdots m_n!}$$ 
so many terms in the vector $\bx^{(m)}$, and hence there are 
${m \choose m_1, m_2, \dots, m_n}!$ permutations for the corresponding columns in $\tilde P$.
Thus, we can generate many new matrices $\tilde P_1$ from $\tilde P$ satisfying
$\tilde P_1 \bx^{(m)} = \bx$ for all $\bx \in \Omega_n$.

\medskip
An interesting question is whether a higher-order Markov chain with transition 
tensor $\tilde P$ satisfying $\tilde P\bx^{(m)} = \bx$ for every $\bx \in \Omega_n$ can be
obtained from the above construction.

\medskip
Next, we turn to higher-order Markov chains satisfying condition (II).

\begin{theorem} Suppose $n > 2$,
$k \in \{1, \dots, n\}$, and $f_k = (e_1 + \cdots + e_k)/k$.
Construct an $m$-th order Markov chain with transition tensor 
$P = (p_{i,i_1,\dots i_m})$ identified as the $n \times n^m$ matrix
$P = [P_1 | \cdots | P_n]$ so that $P_i$ is $n \times n^{m-1}$,
so that every column of $P$ equals $f_k$. Then 
$f_k$ is the only stationary vector of the Markov chain.

\begin{itemize} 
\item[{\rm (1)}] 
If $k = 2$, replace the first volume of  $P_1$ by $e_1$ and all columns of $P_2$ by $e_2$. 
The resulting Markov chain has 2 stationary vectors, namely, 
$e_1$ and $e_2$.

\item[{\rm (2)}]
If $2 < k \le n$, replace the column of $P_i$ indexed by
$(i_1, \dots, i_m) = (i, \dots, i)$
by $e_i$ for all $i = 1, \dots, k$,
and all other columns by
$e_k$ for $i = 1, \dots, k$, then the resulting
Markov chain has $k$ stationary vectors, namely, 
$e_1, \dots, e_k$. 

\item[{\rm (3)}] 
Suppose $k = n$. If we replace the 
the column of $P_i$ indexed by
$(i_1, \dots, i_m) = (i, \dots, i)$ by $e_i$ for all 
$i = 1, \dots, n$, then 
the resulting Markov chain has $n+1$ stationary vectors, namely, 
$e_1, \dots, e_n$ and $f_n$.  

\item[{\rm (4)}] If $k \in \{2, \dots, n-1\}$ and we replace the 
columns of $P_i$ indexed by $(i_1, \dots, i_m)$ with $1 \le i_1, \dots, i_m \le k$
by $e_i$ for $i = 1, \dots, k$, 
then the set of stationary vectors for 
the Markov chain equals
$\conv\{e_1, \dots, e_k\}$.

\item[{\rm (5)}] Suppose $k \in \{2, \dots, n-1\}$ and we reset the matrices
$P_1, \dots, P_n$ so that for each $i = 1, \dots, n$,
the  columns of $P_i$
indexed by $(i_1, \dots, i_m)$ with $1 \le i_1, \dots, i_m \le k$, 
equal 
$$v_i = \begin{cases}  e_i & \hbox{ if } \ $i = 1, \dots, k$,\\
(e_{k+1} + \cdots + e_n)/(n-k) & \hbox{ if } \  $i = k+1, \dots, n$,\end{cases}$$ 
and all other columns of $P_i$ equal $f_n$.
Then the set of stationary vectors for 
the Markov chain equals 
$\{f_n\} \cup \conv \{e_1, \dots, e_k\}$.
\end{itemize}
\end{theorem}

\it Proof. \rm The proof is an easy adaptation of that of Theorem \ref{main2}.\qed

\medskip
To conclude our note, we remark that there are many interesting questions concerning
the stationary vectors of higher-order Markov chains that deserve further study.

\medskip\noindent
{\bf Acknowledgment}

The study of the problem in the paper
began when Li was visiting the University of Hong Kong in the Spring of 2012.
He gracefully acknowledge the support and hospitality of the colleagues at the
Department of Mathematics at the University of Hong Kong. 
The research of Li was supported by USA NSF and HK RCG.
The authors would also like to thank the referee for some 
helpful comments and references.

\end{document}